\documentclass[12pt]{article}

\usepackage{amsmath,amssymb,amsfonts,theorem,makeidx,latexsym,epsfig}
\usepackage{color}
\usepackage{graphicx}
\usepackage{subcaption} 
\usepackage{url}

\newtheorem{defn}{Definition}[section]

\newtheorem{proposition}[defn]{Proposition}

{\theorembodyfont{\rmfamily}

\newtheorem{ex}[defn]{Example}}
\newtheorem{thm}[defn]{Theorem}
\newtheorem{prop}[defn]{Proposition}

\newcommand{\h}{{\cal H}}
\newcommand{\mn}{\mathbb N}
\newcommand{\mr}{\mathbb R}
\newcommand{\mz}{\mathbb Z}
\newcommand{\mq}{\mathbb Q}

\def\dx{{{\rm d}x}}
\def\bp{{\noindent\bf Proof. \ }}
\def\ep{\hfill$\square$\par\bigskip}
\def\vn{\vspace{.1in}\noindent}
\def\<{\langle}
\def\>{\rangle}
\newcommand{\seqgr[1]}{\{#1_k\}_{k=1}^\infty}

\title{On compactly supported dual windows of Gabor frames}

\date{\today}

\author{Diana T. Stoeva  \vspace{.07in} \\ 
{\normalsize Acoustics Research Institute, Austrian Academy of Sciences,} \\
{\normalsize Wohllebengasse 12-14,  Vienna 1040, Austria}\\
{\normalsize dstoeva@kfs.oeaw.ac.at}
}

\begin{document}

\maketitle

\begin{abstract} 
The main purpose of the paper is to give a characterization of all compactly supported dual windows of a Gabor frame. As an application, we consider an iterative procedure for approximation of the canonical dual window via compactly supported dual windows on every step. In particular, the procedure allows to have approximation of the canonical dual window via dual windows from certain modulation spaces or from the Schwartz space. 
\end{abstract}


\vspace{.1in}

\section{Introduction and main results}

Throughout the paper, $\h$ denotes a separable Hilbert space. 
A sequence $G=\seqgr[g]$ with elements from $\h$ is a \emph{frame for $\h$} \cite{DSframe}, if there exist positive constants $A_G$ and $B_G$ such that 
$$A_G\|h\|^2\leq \sum_{k=1}^\infty |\<h,g_k\>|^2 \leq B_G\|h\|^2,  \forall h\in\h.$$
Frames extend orthonormal bases, allowing redundancy, and still guarantee perfect and stable reconstruction. 
Given a frame $G$ for $\h$, there always exists a frame $F=\seqgr[f]$ for $\h$ so that
$$h= \sum_{k=1}^\infty \<h,g_k\> f_k = \sum_{k=1}^\infty \<h,f_k\> g_k, \forall h\in\h;$$
such a frame $F$ is called a {\it dual frame of $G$}. The so called {\it frame operator of $G$} is given by $S_Gh=\sum_{k=1}^\infty \<h,g_k\>g_k$; it is a bounded bijective operator on $\h$ and the sequence $\{S_G^{-1}g_k\}_{k=1}^\infty$ is a dual frame of the frame $G$, called the {\it canonical dual of $G$}. 
  When $G$ is a frame which is not a Schauder basis (so called \emph{overcomplete  
  frame}), it has other dual frames (infinitely many) in addition to the canonical one.

For certain engineering fields, e.g. in signal processing, frames of Gabor structure play an essential role. Given $g\in L^2(\mr)$ and positive constants $a$ and $b$, a \emph{Gabor system} is a system of the form $\{E_{mb}T_{na} g\}_{m,n\in\mz}$; $g$ is called the \emph{window} of the system $\{E_{mb}T_{na} g\}_{m,n\in\mz}$. Here $E_\nu$ and $T_t$ denote the standard modulation\footnote{For $\nu\in\mr$, the modulation operator $E_\nu :L^2(\mr)\to L^2(\mr)$ is determined by $E_\nu g \, (x) =e^{2\pi {\rm i}\nu x} g(x)$, $x\in\mr$.}   and translation\footnote{For $t\in\mr$, the translation operator $T_t:L^2(\mr)\to L^2(\mr)$ is determined by $T_t g \, (x) = g(x-t)$, $x\in\mr$.} operator, respectively.   A \emph{Gabor frame} is a Gabor system which is a frame for $L^2(\mr)$. 
If a Gabor system $\{E_{mb}T_{na} \phi\}_{m,n\in\mz}$ is a dual frame of a Gabor frame 
$\{E_{mb}T_{na} g\}_{m,n\in\mz}$, $\phi$ is called a {\it dual window of} $\{E_{mb}T_{na} g\}_{m,n\in\mz}$. 
  For more on general frame theory, as well as on Gabor frames and other structured frames, we refer e.g. to the books \cite{CasazzaK,Olebook,FS1,FS2,Gbook,Hbook}. 
  
While in general the structure of a given frame is not inherited by the canonical dual (e.g. this is the  case with wavelet frames), 
the canonical dual of a Gabor frame $G=\{E_{mb}T_{na} g\}_{m,n\in\mz}$ is also with Gabor structure - the function $S_G^{-1}g$ is a dual window of $G$ (see, e.g., \cite{Daubechies} 
or \cite[Theorem 12.3.2]{Olebook}) called the {\it canonical dual window}. 
However, the canonical dual of a Gabor frame 
may fail some other nice properties desired for applications (e.g. compact support, smoothness,  
good time-frequency localization, 
 and other). In particular, a Gabor frame, which is also a Schauder basis, and its only dual (the canonical one) can never be well localized in both time and frequency, cf. the Balian-Low Theorem \cite{Daubechies,BHW}. 
In contrast to the above, in the overcomplete case there exist Gabor frames with a nice generating window and nice dual 
windows with respect to certain desired properties (see, e.g., \cite{CGabor,CK,CKK2013,CKK2015,LLM,PHLB}). 
Although not all the dual frames of an overcomplete Gabor have to be with a Gabor structure \cite{Li}, there are infinitely many dual Gabor frames.  
 There exist characterizations of all the dual windows of a Gabor frame via the 
Wexler-Raz biorthogonality relation 
\cite{WR,Jannsen,DLL,FK} 
and via more constructive 
approaches 
 \cite{Li,C2006,HLS} involving the inverse of the frame operator.  
For computational purposes, those of the dual Gabor frames which 
have compactly supported windows  
are of significant importance. 
The main purpose of this paper is to give a characterization of all the 
dual windows which have compact support and to avoid operator inversions: 

\begin{thm}\label{thmcompact}
  Let  $g\in L^2(\mr)$ 
   be compactly supported 
   and such that $G=\{E_{mb}T_{na}g\}_{m,n\in\mz}$ is a frame for $L^2(\mr)$ for some $a,b>0$. 
 Assume that there is a compactly supported function $g^d\in L^2(\mr)$ so that  $G^d=\{E_{mb}T_{na}g^d\}_{m,n\in\mz}$ is a dual frame of $G$. 
Then there is a finite set $K$ (dependent only on $g$, $g^d$, and $a$) so that
 all the compactly supported dual windows of $G$ are  
 in the form 
 \begin{equation}\label{dualallgabcomp}
\phi = g^d+w -\sum_{k\in K}\sum_{j\in\mz} \<g^d, E_{jb}T_{ka} g\> E_{jb}T_{ka}w
\end{equation}
where $w\in L^2(\mr) $ is compactly supported 
and such that $\{E_{mb}T_{na}w\}_{m,n\in\mz}$ is a Bessel sequence in $L^2(\mr)$. 
\end{thm}

For classes of Gabor frames for which the assumptions of the above theorem are fulfilled, see, e.g., \cite{CGabor,CK,CKK,CKK2013,CKK2015}.

\vspace{.1in}
As an application of Theorem \ref{thmcompact}, we can give a procedure for approximation of the canonical dual window via compactly supported dual windows (Proposition \ref{propgaborcomp}). 
Recall that the canonical dual window of a Gabor frame  $\{E_{mb}T_{na}g\}_{m,n\in\mz}$  has some nice properties, 
 e.g. it has minimal $L^2$-norm among the dual windows \cite[Prop. 3.2]{JanssenCanDualWindow}; it belongs to the modulation space $M_{1,1}^v$ when 
 $g\in M_{1,1}^v$, $ab\in\mq$, and $v$ is a polynomial weight \cite[Cor. 3.5]{FG}; 
it belongs to the Schwartz space $\mathcal{S}$ when 
$g\in\mathcal{S}$ 
\cite[Prop. 5.5]{JanssenCanDualWindow} 
 \cite[Cor. 3.6]{FG}. 
 However, often it is difficult to find the canonical dual in explicit form and for such cases  it is of interest to consider approximations via dual frames which are explicitly determined. 
Furthermore, the canonical dual window of a Gabor frame with a compactly supported window is not necessarily with a compact support  
(though there are Gabor frames for which the canonical dual  is compactly supported, see e.g. \cite{BW,Boe}, \cite[Ex. 12.3.3]{Olebook}). 
For such Gabor frames it is of interest to have approximations of the canonical dual via compactly supported dual windows. 

Note that the known frame algorithm and its accelerations based on the Chebyshev method and the conjugate gradient method \cite{GrAlg} 
for approximation of the canonical dual of a general frame, as well as the variant of the Schulz iterative algorithm \cite[Alg. IV]{JanssenSchulzAlg} 
for approximation of the canonical dual window of a Gabor frame, 
 and various finite section methods \cite{Strohmer,CC}, 
use iteration steps which are not necessarily dual frames. 
As far as we are aware, the only paper which concerns an approximation of the canonical dual window via exact dual windows 
is \cite{Tobias}. 
The algorithm in \cite{Tobias} is designed for Gabor frames with very specific windows (totally positive functions and exponential B-splines) and it provides compactly supported dual windows which converge exponentially 
 to the canonical dual window. 
 Here we consider 
 an algorithm (variation of the frame algorithm), 
  which provides compactly supported dual windows on every iteration and applies to a 
  quite general class of Gabor frames.

\begin{prop} \label{propgaborcomp} 
Let $g, g^d\in L^2(\mr)$ and $a,b>0$ be such that 
 $G=\{E_{mb}T_{na}g\}_{m,n\in\mz}$ is a frame for $L^2(\mr)$ and  $G^d=\{E_{mb}T_{na}g^d\}_{m,n\in\mz}$ is a dual frame of $G$. 
 Let $\lambda$ be a positive constant such that $\|I-\lambda S_G\|<1$. Consider the sequence $\{g^p\}_{p=0}^\infty$ determined as follows:
\begin{equation}\label{dualngab}
\begin{aligned}
g^0&:=  g^d,\\
 g^{p+1}&:=  \lambda g + (I - \lambda S_G) g^p, \ p\in\mn_0. 
\end{aligned}
\end{equation}
Then $\{g^p\}_{p=0}^\infty$ is a sequence of dual windows of $G$ which converge to the canonical dual window  $S_G^{-1}g$. 
If $g$ and $g^d$ are furthermore compactly supported, then $g^p$, $p\in\mn$, are also compactly supported. 
\end{prop}

\noindent
{\bf Remark.} Note that under the assumptions of the above proposition, if $ab\in\mq$ and $g$ and $g^d$ belong to the modulation space $M_{1,1}^{\rm w}(\mr)$ \cite{Feicht83} for a 
 submultiplicative weight ${\rm w}$ of polynomial growth 
 (resp. to the Schwartz space $\mathcal{S}(\mr)$), then applying \cite[Theorem 3.4  
 and (S6)]{FG} and using induction on $p$ we get that the functions $g^p$, $p\in\mn$, given by (\ref{dualngab}), belong to $M_{1,1}^{\rm w}(\mr)$ (resp. to $\mathcal{S}(\mr)$).

\vspace{.1in}
The paper is organized as follows. 
In Section 2 we provide the necessary background results needed to prove Theorem \ref{thmcompact}.  
These are characterizations of all the dual frames based on any given dual frame 
and they can be of independent interest  
 serving as a tool for construction of 
variety of dual frames. 
This can be very useful 
e.g. in cases when one 
 searches for dual frames satisfying 
 additional constraints motivated by applications. 
Section 3 is devoted to 
 application of 
  Section 2 for approximation of the canonical dual frame. 
  A proof of Proposition \ref{propgaborcomp} is given and a procedure applying to general frames 
   is also considered. 
In both Sections 2 and 3 we give illustrative examples and provide implementations under the Matlab environment.

\section{Characterizations of dual frames}\label{sec2}
The main focus of this paper is Theorem \ref{thmcompact} - a characterization of all the compactly supported dual windows of a Gabor frame. In order to prove this theorem, we will 
need the following results (Propositions \ref{prop1} and \ref{dualgabfr}), 
which are of independent interest as well. 
As a first step, 
observe 
that the known characterization of all the dual frames via the canonical dual 
 \cite{Li,C2006,HLS} 
can be extended to the use of any dual frame instead of the canonical dual.

\begin{proposition}\label{prop1}
Let $G=\{g_k\}_{k=1}^\infty$ be a frame for $\h$ and let $G^d=\{g_k^d\}_{k=1}^\infty$ be a dual frame of $G$. 
Then all the dual frames of $G$ are precisely the sequences
\begin{equation}\label{dualall}
\{ g_k^d + w_k-\sum_{j=1}^\infty \<g_k^d, g_j\> w_j \}_{k=1}^\infty,
\end{equation}
where $W=\{w_k\}_{k=1}^\infty$ is a Bessel sequence in $\h$.
\end{proposition}
\bp A proof can be done in a similar way as the proof  
of characterizations based on the canonical dual (see, e.g., \cite[Theor. 6.3.7]{Olebook}). For   
convenience of the readers, let us give a sketch of the proof. 
First take a sequence in the form (\ref{dualall}) with  $W=\{w_k\}_{k=1}^\infty$ being a Bessel sequence in $\h$ and write it in the form 
\begin{equation*}\label{op1}
\{(T_{G^d}+T_W-T_WU_GT_{G^d})\delta_k\}_{k=1}^\infty,
\end{equation*}
where $\delta_k$ means the $k$-th canonical vector, $k\in\mn$. 
Denote $V:=T_{G^d}+T_W-T_WU_GT_{G^d}$. Clearly, $V$ is a bounded operator from $\ell^2$ into $\h$ and $VU_G= Id_\h$. 
Now use the known result that 
the dual frames of $G$ are precisely the sequences of the form  $\{L\delta_k\}_{k=1}^\infty$, where $L:\ell^2\to\h$ is a bounded left inverse of $U_G$ (see, e.g., \cite[Lemma 6.3.5]{Olebook}) and conclude that  $\{V\delta_k\}_{k=1}^\infty$ is a dual frame of $G$. 

Conversely, take a dual frame $F=\{f_k\}_{k=1}^\infty$ of $G$. The sequence  $F$ can be written in the form (\ref{dualall})  using for example $W=F$.
\ep

In the case of Gabor frames, we can also use the known characterization of the  dual Gabor frames based on the canonical dual \cite{Li,C2006,HLS} and extend it to the use of any other dual Gabor frame instead of the canonical one. 
This can be useful for frames for which an explicit formula for the canonical dual window is hard to get, but another dual window is known in explicit form (for such cases see, e.g., \cite{CKK2015},\cite[Sec. 12.5,12.6]{Olebook}). It is also useful for the characterization of all compactly supported dual windows (Theorem \ref{thmcompact}), which could not have been done using only the known characterization based on the canonical dual, because the canonical dual window of a Gabor frame with a compactly supported generating window is not necessarily compactly supported.

\begin{proposition} \label{dualgabfr} 
Let $g, g^d\in L^2(\mr)$ and $a,b>0$ be such that 
 $G=\{E_{mb}T_{na}g\}_{m,n\in\mz}$ is a frame for $L^2(\mr)$ and  $G^d=\{E_{mb}T_{na}g^d\}_{m,n\in\mz}$ is a dual frame of $G$. 
Then all the dual windows of $G$ are the functions of the form
\begin{equation}\label{dualallgab}
\phi=g^d+w -\sum_{j,k\in\mz} \<g^d, E_{jb}T_{ka} g\> E_{jb}T_{ka}w,
\end{equation}
where $w\in L^2(\mr) $ is such that $\{E_{mb}T_{na}w\}_{m,n\in\mz}$ is a Bessel sequence in $L^2(\mr)$. 
\end{proposition}
\bp
A proof can be done in a similar way as the proof of \cite[Prop. 12.3.6]{Olebook}, using Proposition \ref{prop1} instead of \cite[Theorem 6.3.7]{Olebook}. We give a sketch of the proof for convenience of the readers. 
 
For one of the directions, take $w\in L^2(\mr) $  such that $\{E_{mb}T_{na}w\}_{m,n\in\mz}$ is a Bessel sequence in $L^2(\mr)$. Then the  function $\phi$ determined by (\ref{dualallgab}) is well defined and belongs to $L^2(\mr) $. 
Using similar calculations as in \cite[Lemma 12.3.1]{Olebook}, the Gabor system $\{ E_{mb}T_{na}\phi\}_{m,n\in\mz}$  can be written as
$$
\{ E_{mb}T_{na}g^d+E_{mb}T_{na} w -  \sum_{j,k\in\mz} \<E_{mb}T_{na}g^d, E_{jb}T_{ka} g\> E_{jb}T_{ka}w) \}_{m,n\in\mz},
$$ 
which is a dual frame of $G$ by Proposition \ref{prop1}. 

Conversely, if $\phi$ is a dual window of $G$,  then it can be written in  the form (\ref{dualallgab})  taking e.g. $w=\phi$. 
\ep

We can now prove the characterization of all compactly supported dual windows of a Gabor frame (in the cases when such ones exist):

\vn
{\bf Proof of Theorem \ref{thmcompact}.} First observe that since $g$ and $g^d$   are compactly supported, 
one can determine a finite set $K$ (dependent on $g$, $g^d$, and $a$) so that 
 $ \<g^d, E_{jb}T_{ka} g\>\neq 0$ if and only if $k\in K$. 

Take  a compactly supported $w\in L^2(\mr) $ such that $\{E_{mb}T_{na}w\}_{m,n\in\mz}$ is a Bessel sequence in $L^2(\mr)$. 
 Consider the $L^2$-function $\phi$ given by  (\ref{dualallgabcomp}). 
Clearly,  $\phi$ is also with a compact support. 
   Since the sum on $k$  
 in (\ref{dualallgab}) is actually a sum on $k\in K$, $\phi$ can be written as in 
 (\ref{dualallgab}). Thus, by Proposition \ref{dualgabfr}, $\{E_{mb}T_{na}\phi\}_{m,n\in\mz}$ is a dual frame of $G$. 

Conversely, if $\phi$ is a compactly supported dual window $G$,  then $\phi$ can be written in  the form (\ref{dualallgabcomp})  using for example $w=\phi$. 
\ep

For implementation purposes, it is convenient to use real-valued windows. 
The following statement gives a characterization of all the real-valued compactly supported dual windows. 

\begin{prop}\label{thmcompact2} 
Under the assumptions of Theorem \ref{thmcompact}, let $g$ and $g^d$ be furthermore real-valued. 
For $j,k\in\mz$, denote $P_{j,k}(f) (x):=(\int g^d(x) f(2\pi j bx) T_{ka} g(x)\dx) f(2\pi j bx).$ 
Then there is a finite set $K$ (dependent only on $g$, $g^d$, and $a$) so that
 all real-valued compactly supported dual windows of $G$ are 
in the form 
  \begin{equation}\label{dualallgab3}
\phi = g^d+w -\sum_{k\in K} 
(\<g^d, T_{ka} g\>+2\sum_{j\in\mn} 
(P_{j,k}(cos) + P_{j,k}(sin))) T_{ka}w
\end{equation}
where $w\in L^2(\mr) $ is real-valued compactly supported 
and such that $\{E_{mb}T_{na}w\}_{m,n\in\mz}$ is a Bessel sequence in $L^2(\mr)$. 
\end{prop}
\bp  
 As in Theorem \ref{thmcompact}, there is a finite set $K$ (dependent on $g$, $g^d$, and $a$) so that  
 $ \<g^d, E_{jb}T_{ka} g\>\neq 0$ if and only if $k\in K$.  
  Take a real-valued compactly supported 
    $w\in L^2(\mr) $ such that $\{E_{mb}T_{na}w\}_{m,n\in\mz}$ is a Bessel sequence in $L^2(\mr)$.  
Using the unconditional convergence of
$\sum_{j,k\in\mz} \<g^d, E_{jb}T_{ka} g\> E_{jb}T_{ka}w$, one can proceed to obtain  representation with real-valued functions: 
\begin{eqnarray*} 
 & & \sum_{k\in K} \sum_{j\in\mz} \<g^d, E_{jb}T_{ka} g\> E_{jb}T_{ka}w \\ 
&=& \sum_{k\in K} \left(\<g^d, T_{ka} g\> +
\sum_{j\in\mn}(\<g^d, E_{jb}T_{ka} g\> E_{jb}
+ 
\<g^d, E_{-jb}T_{ka} g\> E_{-jb})\right)T_{ka}w.
\\
&=& 
\sum_{k\in K} \left(\<g^d, T_{ka} g\> +
2\sum_{j\in\mn}(P_{j,k}(cos) + P_{j,k}(sin))\right)T_{ka}w.
\end{eqnarray*} 
Now let $\phi$ be 
determined by (\ref{dualallgab3}) and thus being real-valued. By the above, it follows that $\phi$ satisfies (\ref{dualallgabcomp}) and hence, by Theorem \ref{thmcompact}, $\phi$ is a compactly supported dual window of $G$. 

For the converse part, if $\phi$ is a real-valued compactly supported dual window of $G$, then taking $w=\phi$ one can write $\phi$ as in (\ref{dualallgabcomp})  and hence as in   (\ref{dualallgab3}). 
\ep

\begin{figure}[h!]
  \centering
  \begin{subfigure}[b]{0.45\linewidth}
    \includegraphics[width=\linewidth]{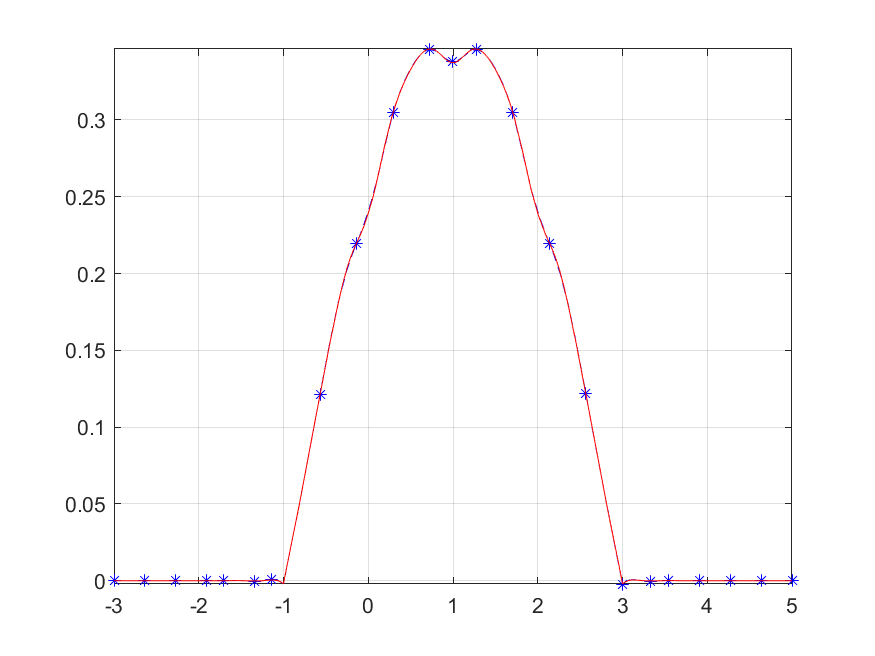}
    \caption{ using $\lambda=\frac{1}{10}$}
  \end{subfigure}
  \begin{subfigure}[b]{0.45\linewidth}
    \includegraphics[width=\linewidth]{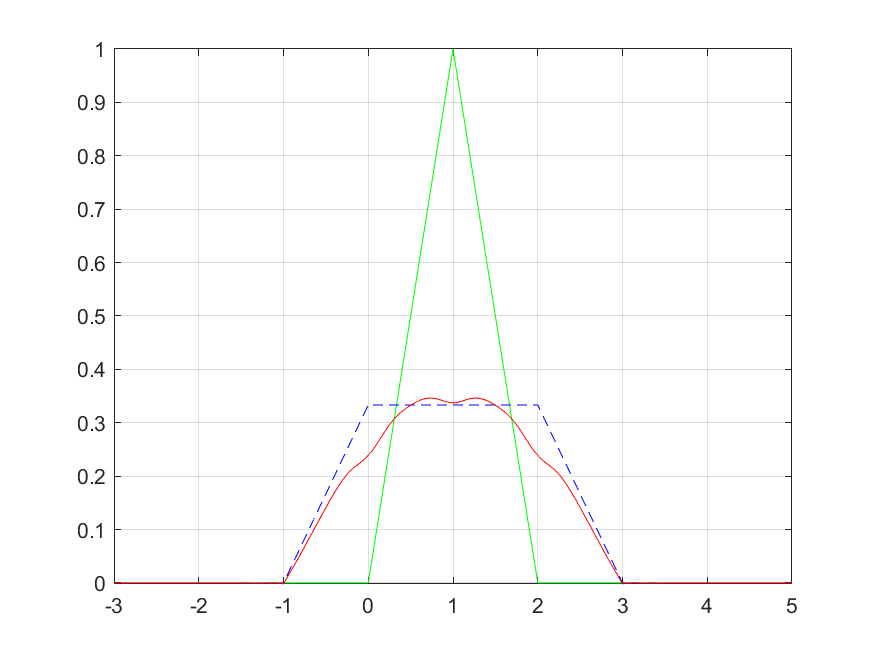}
    \caption{ using $\lambda=\frac{1}{10}$}
  \end{subfigure}
   \begin{subfigure}[b]{0.45\linewidth}
    \includegraphics[width=\linewidth]{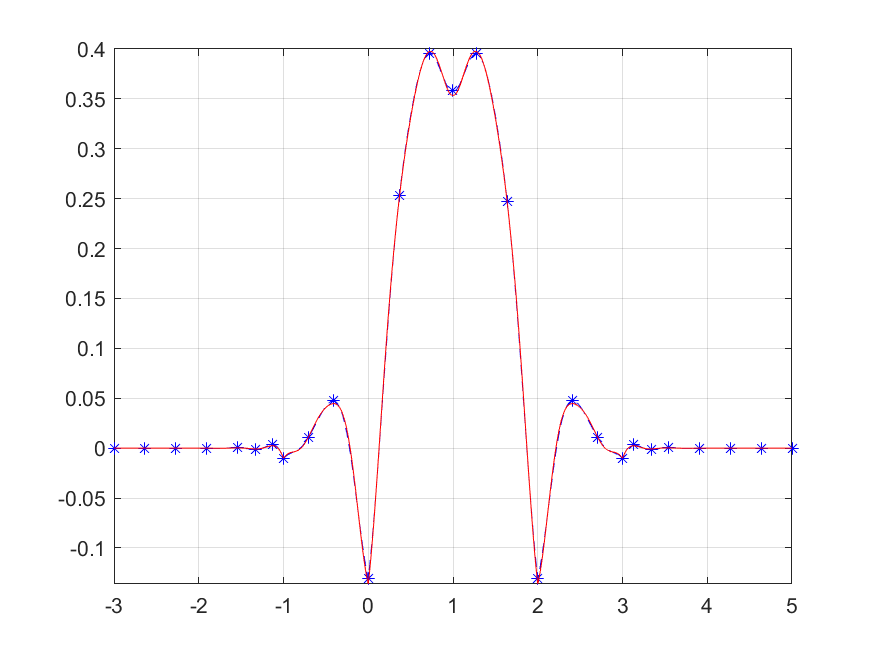}
    \caption{ using $\lambda=\frac{1}{2}$}
  \end{subfigure}
     \begin{subfigure}[b]{0.45\linewidth}
    \includegraphics[width=\linewidth]{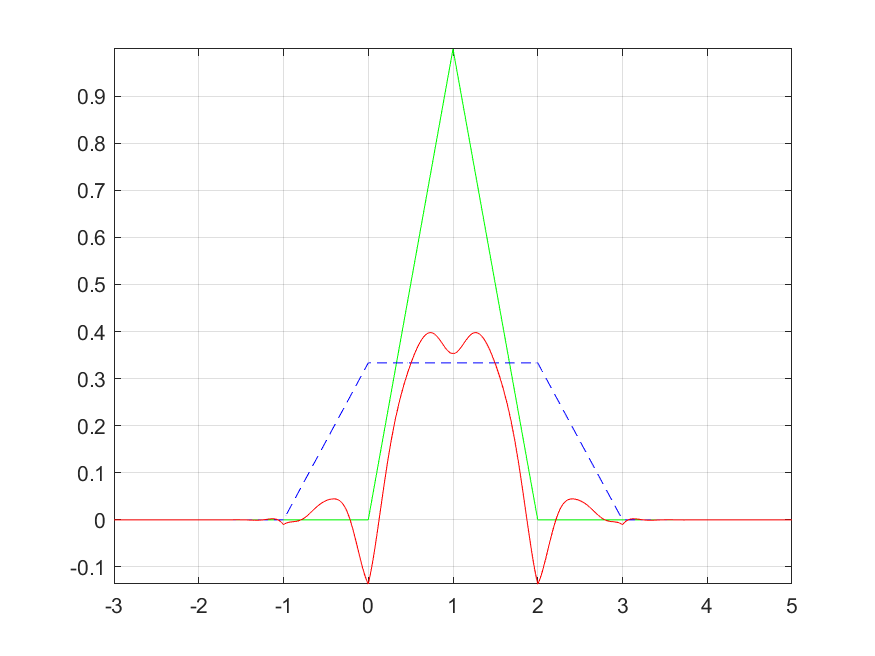}
    \caption{ using $\lambda=\frac{1}{2}$}
  \end{subfigure}
  \caption{Visualization of $g$, $g^d$, and 
   $\phi$ from Example \ref{exb2tilde}. For (a) and (c):  
  comparison of $\phi$ calculated with sum on $j$ truncated to $6$ terms (in blue with stars) and to $7$ terms (in red). For (b) and (d):  the function $g$   (in green), $g^d$ (in blue), and $\phi$ calculated with sum on $j$ truncated to $7$ terms (in red).}
  \label{figGab}
\end{figure}

As an illustration of Theorem \ref{thmcompact} and Proposition \ref{thmcompact2}, consider Example \ref{exb2tilde} and Fig.\,\ref{figGab}. 
The Matlab script which was used to produce Fig.\,\ref{figGab} can be found at \url{https://www.oeaw.ac.at/isf/ondualframes}. 

\begin{ex} \label{exb2tilde}
Let $a=1$ and $b=1/3$. Consider the Gabor frame $G=\{E_{mb} T_{na}g\}_{m,n\in\mz}$, where $g$ is the B-spline 
$$B_2(x)=\left\{
\begin{array}{ll}
x, & x\in [0,1),\\
2-x, & x\in [1,2],\\
0, & x\notin [0,2],
\end{array}
 \right.
   $$
and the function
$$h_2(x)=\left\{
\begin{array}{ll}
\frac{1}{3}(x+1), & x\in [-1,0),\\
\frac{1}{3}, &x\in [0,2),\\
1-\frac{1}{3}x, & x\in [2,3],\\
0, & x\notin [-1,3],
\end{array}
 \right.
   $$
which is a dual window of $G$ \cite{CK}. 
Consider the function $\phi$ given by 
 (\ref{dualallgabcomp}) 
with $g^d=h_2$ and $\omega=\lambda B_2$, where $\lambda$ is a positive constant. 
Observe that the set $K$ in  (\ref{dualallgabcomp}) 
 is $\{-2,-1,0,1,2\}$. 
By Theorem \ref{thmcompact}, 
$\phi$ is a dual window of $G$. In Fig.\,\ref{figGab} we visualize  
$\phi$ for two values of $\lambda$. 
 \end{ex}
 
Concerning the  implementation of 
Example \ref{exb2tilde}, 
note that the visualized $\phi$ on Fig.\,\ref{figGab} 
(calculated using truncated  sum on $j$ in (\ref{dualallgab3}))
 is symmetric with respect to the line $x=1$, like the given $g$ and $g^d$. Actually, 
 in  the setting of Prop. \ref{thmcompact2}, one can show that if $g$, $g^d$, and $w$ are furthermore symmetric with respect to some line $x=x_0$ ($x_0\in\mr$), then any truncated representation of $\phi$ by  (\ref{dualallgab3})  with a finite sum on $j$ gives a symmetric function with respect to $x=x_0$.

\vspace{.1in}
Clearly, 
one can consider the case of frames of translates and to write respective characterizations of the dual frames of translates in the spirit of Proposition \ref{dualgabfr} and Theorem \ref{thmcompact}.

\section{Approximation of the canonical dual frame with dual frames on every step}

As motivated in the introduction, our interest in this section is in approximation of the canonical dual frame via dual frames. 
Let us begin with the Gabor case and the approximation of the canonical dual window via compactly supported dual windows.

\vn
{\bf Proof of Proposition \ref{propgaborcomp}: }
Consider the sequence $\{g^p\}_{p=0}^\infty$ given by (\ref{dualngab}). 
For a given 
$w\in L^2(\mr) $ such that $\{E_{mb}T_{na}w\}_{m,n\in\mz}$ is a Bessel sequence in $L^2(\mr)$, consider the sequence $q^p$, $p\in\mn_0$, determined as follows:
$q^0 :=  g^d$ 
and 
$ q^{p+1}:=  q^p+w -\sum_{j,k\in\mz} \<q^p, E_{jb}T_{ka} g\> E_{jb}T_{ka}w$, $p\in\mn_0$.  
Using Proposition \ref{dualgabfr}, it follows by induction that $q^p$ is a dual window of $G$ for every $p\in\mn_0$.   
Choosing $w=\lambda g$, the sequence $\{q^p\}_{p\in\mn_0}$ becomes the same as $\{g^p\}_{p\in\mn_0}$.  
Now the representation 
$$g^p - S_G^{-1}g =  (I-\lambda S_G)^p (g^d- S_G^{-1}g), \ p\in\mn_0,$$ 
 leads to
 $\|g^p - S_G^{-1}g\| \leq \|I-\lambda S_G\|^p (\sqrt{B_{G^d}}+ \frac{1}{\sqrt{A_G}})\to 0\ \mbox{as $p\to\infty$}.$
  
 When $g$ and $g^d$ are furthermore compactly supported, using Theorem \ref{thmcompact} one can conclude that every dual window $g^p$, $p\in\mn$, is compactly supported. 
\ep

\vspace{.1in} 
For general frames, the following procedure holds for approximation of the canonical dual via dual frames:

\begin{prop}\label{approxcandual}
Let $G=\{g_k\}_{k=1}^\infty$ be a frame for $\h$ and let $G^d=\{g^d_k\}_{k=1}^\infty$ be a dual frame of $G$.  Let $\lambda$ be a positive constant such that $\|I-\lambda S_G\|<1$. Consider the sequence $F^p=\{f^p_k\}_{k=1}^\infty, p\in\mn_0,$ determined as follows: 
\begin{eqnarray}
F^0&:=&  G^d, \label{dual0}\\
 F^{p+1}&:=&  \{ \lambda g_k + f^p_k -\lambda S_G f^p_k \}_{k=1}^\infty, \ p\in\mn_0.\label{dualnstep}
\end{eqnarray}
Then  $F^p$ is a dual frame of $G$ for every $p\in\mn_0$ and 
$\lim_{p\to\infty}F^p=\widetilde{G}$ 
uniformly on $k$, i.e.,  
 for every $\varepsilon >0$, there exists $N_\varepsilon\in\mn$ so that $\|f^p_k-\widetilde{g}_k\|<\varepsilon$ for every $p>N_\varepsilon$ and every $k\in\mn$. 
\end{prop}

\bp For a given Bessel sequence $W=\{w_k\}_{k=1}^\infty$ in $\h$, consider the sequence $Q^p=\{q^p_k\}_{k=1}^\infty, n\in\mn_0,$ determined as follows:
$Q^0 :=  G^d$ 
and 
$ Q^{p+1}:=  \{q^p_k + w_k-\sum_{j=1}^\infty \<q^p_k, g_j\> w_j \}_{k=1}^\infty$, $n\in\mn_0$.  
Using Proposition \ref{prop1}, it follows by induction that $Q^p$ is a dual frame of $G$ for every $p\in\mn_0$.  
Choosing $W=\lambda G$, the sequence $\{Q^p\}_{p\in\mn_0}$ becomes $\{F^p\}_{p\in\mn_0}$.  

As in the frame algorithm, for every $p,k\in\mn$ one can write 
$f^p_k - S_G^{-1}g_k 
= (I-\lambda S)^p (f^0_k -S^{-1}g_k).
$
Therefore, 
$$
\|f^{p}_k - S_G^{-1}g_k\|
\leq 
\|I-\lambda S\|^p  (\sqrt{B_{G^d}}+ \frac{1}{\sqrt{A_G}}), \ \forall k,p\in\mn,
$$
which leads to the desired conclusion. 
\ep

Note that the consideration of approximation of the canonical dual via dual frames in this section was naturally motivated by the results in 
Section \ref{sec2}, but it turned out to be simply related to the classical frame algorithm
with difference just in the initial step. 
The initialization with a null vector in the classical frame algorithm (see, e.g., \cite[Alg. 5.1.1]{Gbook})
 actually prevents the sequences from the next steps of the algorithm to be dual frames -  the first step leads to the sequence $\lambda G$, which can be a dual frame of $G$ only if $G$ is a $\frac{1}{\lambda}$-tight frame. 
By Proposition \ref{approxcandual}, 
using a dual frame of $G$ in the initialization step of this algorithm is the key to guarantee dual frames on all steps.

\begin{figure}
  \centering
  \begin{subfigure}[b]{0.45\linewidth} 
    \includegraphics[width=\linewidth]{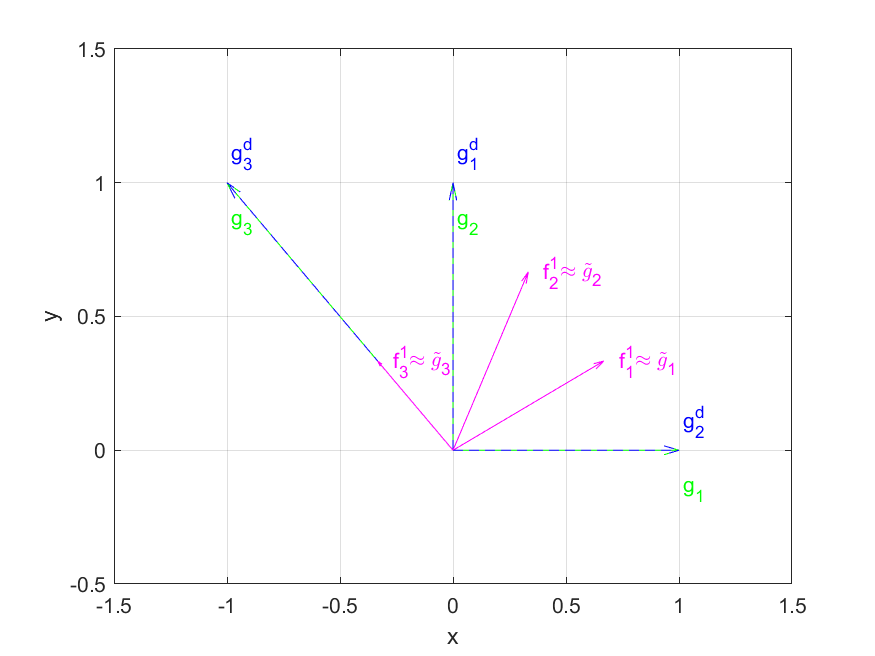}
        \caption{Ex. \ref{ex21}(i) with $\lambda=1/3$}
  \end{subfigure}
  \begin{subfigure}[b]{0.45\linewidth} 
       \includegraphics[width=\linewidth]{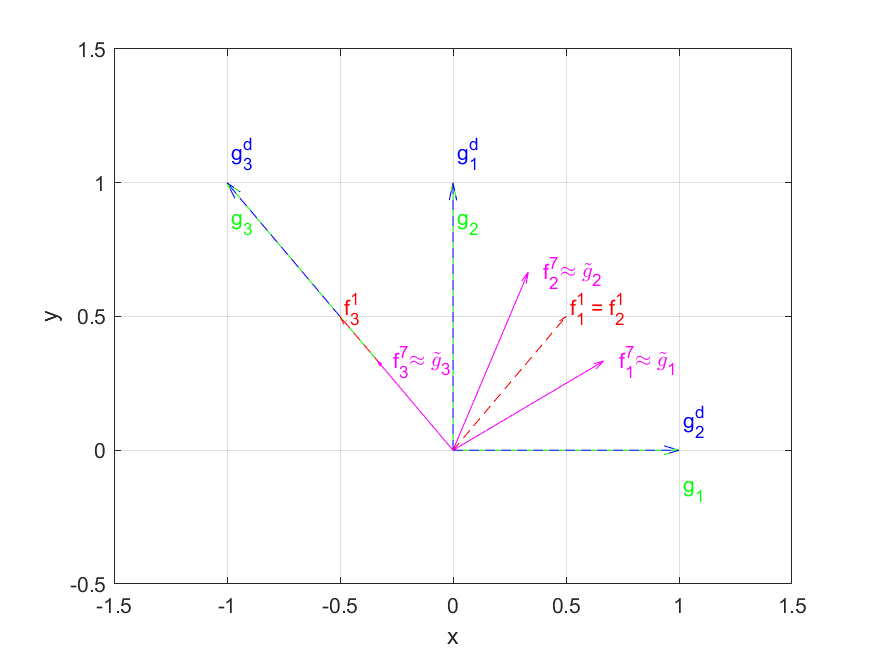}
        \caption{Ex. \ref{ex21}(i) with $\lambda=1/4$}
  \end{subfigure}
    \begin{subfigure}[b]{0.45\linewidth} 
    \includegraphics[width=\linewidth]{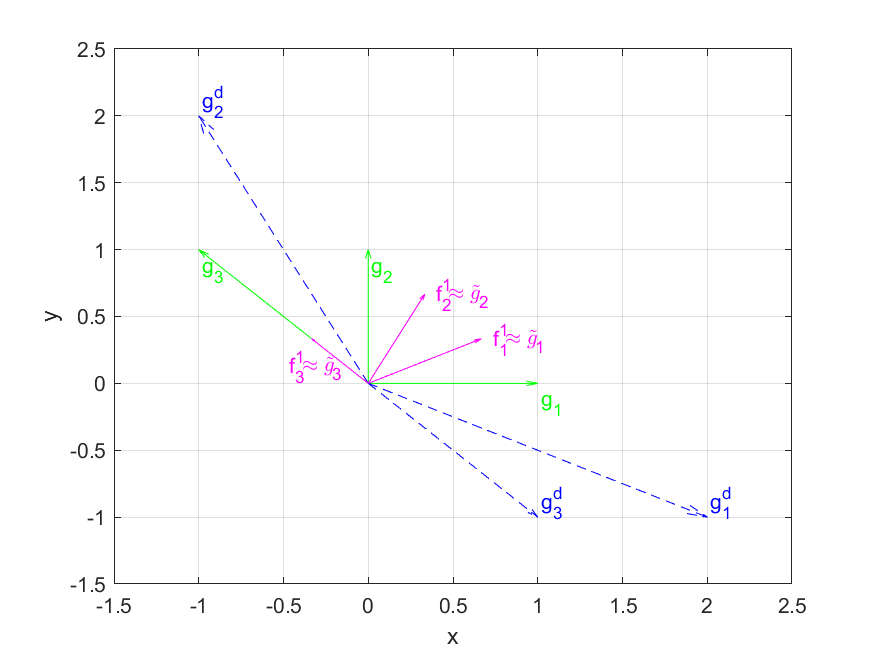}
    \caption{Ex. \ref{ex21}(ii) with $\lambda=1/3$}
  \end{subfigure}
    \begin{subfigure}[b]{0.45\linewidth} 
  \includegraphics[width=\linewidth]{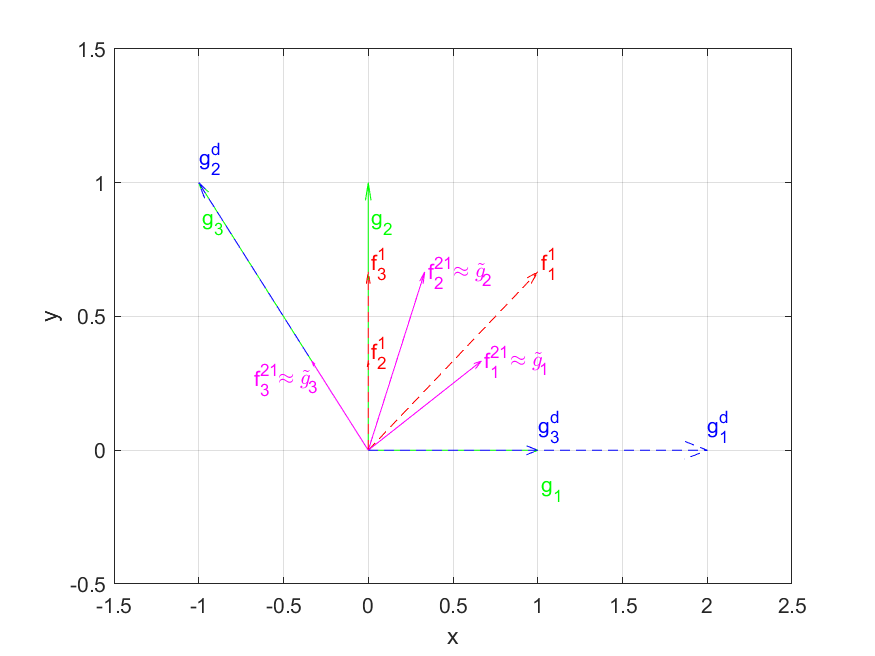}
    \caption{Ex. \ref{ex21}(iii) with $\lambda=1/3$}
  \end{subfigure}
    \begin{subfigure}[b]{0.45\linewidth}
  \includegraphics[width=\linewidth]{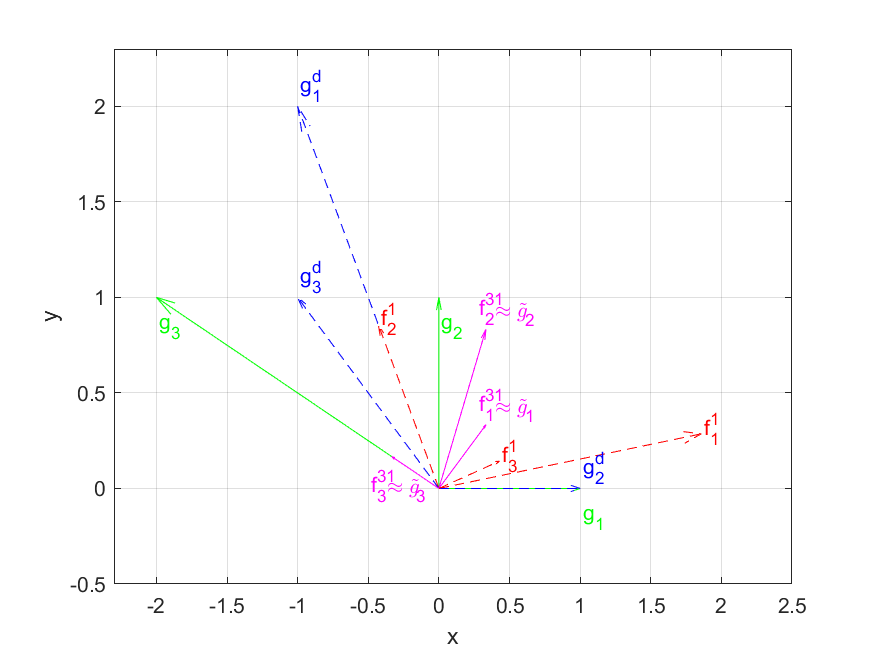}
    \caption{Ex. \ref{ex21}(iv) with $\lambda=2/7$}
  \end{subfigure}
    \begin{subfigure}[b]{0.45\linewidth} 
  \includegraphics[width=\linewidth]{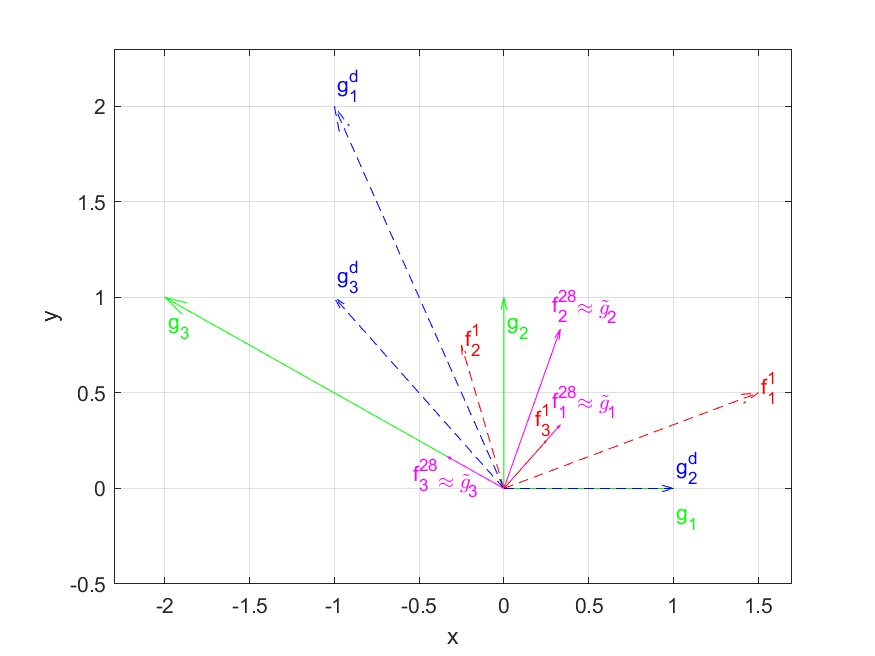}
    \caption{Ex. \ref{ex21}(iv) with $\lambda=1/4$}
  \end{subfigure}
  \caption{The frames from Example \ref{ex21}. In green - the given frame $G$, in blue - the given dual frame $G^d$, in red - the new dual frame from the first iteration step, in magenta - the approximation of the canonical dual $\widetilde{G}$ up to 4 digits after the dot (the number $p$ of the iterations is noted as an upper index of $f$).}
  \label{figapprcandual}
\end{figure}

\vspace{.1in} 
As a simple illustration of Proposition \ref{approxcandual}, consider Example \ref{ex21} and Fig.\,\ref{figapprcandual}. The Matlab scripts for the implementation of Proposition \ref{approxcandual} and for producing Fig.\,\ref{figapprcandual} are available at \url{https://www.oeaw.ac.at/isf/ondualframes}.

\begin{ex}\label{ex21}  
Consider the two-dimensional Euclidean space $\mr^2$, a frame $G$ for $\mr^2$ 
 and its dual frame $G^d$ as given below. 
Based on Proposition \ref{approxcandual}, we 
consider an approximation of the canonical dual frame with precision of the first 4 digits after the decimal dot. The letter $p$ indicates the number of iterations the algorithm took to finish with the desired precision. 
\begin{itemize}

\vspace{-.07in}
\item[{\rm (i)}]  
Consider $G=(e_1, e_2, e_2-e_1)$ and  $G^d=(e_2, e_1, e_2-e_1)$. 
Running the algorithm with 
$\lambda$ being 1/3, 1/4, and $ 1/2(=2/(A_G^{opt}+B_G^{opt})$), 
 the respective values of $p$ are 1, 7, and 14.

\vspace{-.07in}
\item[{\rm (ii)}]  
Consider $G$ from (i) 
and $G^d=(2e_1-e_2, 2e_2-e_1, e_1-e_2)$. 
Running the algorithm with 
$\lambda=1/3, 1/4, 1/2$, we get $p=1, 8,15$, respectively.

\vspace{-.07in}
\item[{\rm (iii)}]  
Consider $G$  from (i) 
and $G^d=(2e_1, e_2-e_1, e_1)$. 
Running the algorithm with 
$\lambda=1/3, 1/4, 1/2$,  we get $p= 21, 29, 15$, respectively.

\vspace{-.07in}
\item[{\rm (iv)}]  
Consider 
$G=(e_1, e_2, e_2-2e_1)$ and 
 $G^d=(2e_2-e_1, e_1, e_2-e_1)$. 
Running the algorithm with 
$\lambda$ being 1/4, 2/9, and $ 2/7(=2/(A_G^{opt}+B_G^{opt})$), the respective values of $p$ are 28, 32, and 31. 
\ep
\end{itemize}
\end{ex}

\vspace{-.15in}
As the above examples show and as it is natural to expect, the efficiency of the algorithms in Propositions \ref{propgaborcomp} and \ref{approxcandual} depends much on the initial dual frame $G^d$ and $\lambda$. The value $2/(A^{opt}_G+B^{opt}_G)$, which is 
an optimal one for $\lambda$  for some algorithms,  
might be much less efficient in the present algorithm compare to other values of $\lambda$, see Example \ref{ex21}(i)(ii). 
As a brief comparison to the classical frame algorithm, if we run
(\ref{dualnstep}) with initialization $F^0={0}$ 
for the frame $G$ from Example  \ref{ex21}(i) 
 with $\lambda=1/3, 1/4, 1/2$, then the respective values of $p$ are $21, 29, 14$ - compare to the values of $p$ in Example \ref{ex21}(i)(ii). 
It will be the purpose of further work to investigate deeper the efficiency of the present algorithm in dependence of  $G^d$ and $\lambda$, especially in the case of Gabor frames.

\newpage
\noindent 
{\bf Acknowledgements} 
The author is grateful to Ole Christensen for 
valuable ideas and suggestions 
 on the topic. 
She is also thankful for the hospitality of the Technical University of Denmark, where the work on the present topic began. The author also acknowledges support 
from the Austrian Science Fund (FWF) through the START-project
``FLAME" Y 551-N13 and from the Vienna Science
and Technology Fund (WWTF) through Project VRG12-009.


\begin{thebibliography}{10}

\frenchspacing


 

 \bibitem{BHW} J. J. Benedetto, C. Heil, D. F. Walnut: \textit{Gabor systems and the Balian-Low Theorem}. In 
H. G. Feichtinger and T. Strohmer (Eds.), Gabor Analysis and Algorithms. Theory and Applications. Boston, MA: Birkhäuser. Applied and Numerical Harmonic Analysis. 85--122 (1998).



 
 \bibitem{BW} J. J. Benedetto, D. F. Walnut: \textit{Gabor frames for $L^2$ and related spaces}. In: J. J. Benedetto and M. W. Frazier (Eds.), Wavelets: Mathematics and Applications, 
   CRC Press, Boca Raton, FL, 97--162 (1994).
 
 

 
 \bibitem{Boe} H. B\"{o}lcskei, \textit{A necessary and sufficient
condition for dual Weyl-Heisenberg frames to be compactly supported}.  
J. Fourier Anal. Appl. 5(5), 409--419 (1999).

\bibitem{CasazzaK}
P. G. Casazza, G. Kutyniok (Eds.): \textit{Finite Frames. Theory and Applications.} 
Birkh\"auser, Basel (2013).



\bibitem{CC}
P. G. Casazza, O. Christensen: \textit{Approximation of the inverse frame operator and applications to Gabor
frames}. J. Approx. Theory 103(2), 338--356 (2000).


\bibitem{C2006}
O. Christensen: \textit{Frames and generalized shift-invariant systems}. In: P. Boggiatto, L. Rodino, J. Toft, M. W. Wong (Eds.), Pseudo-Differential Operators and Related Topics. Operator Theory: Advances and Applications, vol 164, pp. 193--209, Birkh\"auser, Basel (2006).

\bibitem{CGabor} O. Christensen: \textit{Pairs of dual Gabor frames with compact support
and desired frequency localization}. Appl. Comput. Harmon. Anal. 20, 403--410 (2006).

\bibitem{Olebook}
O. Christensen:   \textit{An Introduction to Frames and Riesz Bases.} Second Expanded Edition, 
 Series: Applied and Numerical Harmonic Analysis, 
Birkh\"auser, Boston (2016).



\bibitem{CK} O. Christensen, R.Y.  Kim: \textit{On dual Gabor frame pairs generated by
polynomials}. J. Fourier Anal. Appl. 16, 1--16 (2010).

\bibitem{CKK} O. Christensen, H. O. Kim, R. Y. Kim: \textit{Gabor windows supported on $[-1,1]$ and compactly supported
dual windows}.  
Appl. Comput. Harmon. Anal. 28, 89--103 (2010).

\bibitem{CKK2013} O. Christensen, H. O. Kim, R. Y. Kim: \textit{Regularity of dual Gabor windows}.  Abstr. Appl. Anal., vol. 2013, Article ID 747268 (2013).



\bibitem{CKK2015} O. Christensen,  H.O. Kim, R.Y. Kim: \textit{On entire functions restricted
to intervals, partition of unities, and dual Gabor frames}.
Appl. Comput.Harmon.Anal. 38, 72--86 (2015).


\bibitem{Daubechies}
I. Daubechies: \textit{The wavelet transform, time-frequency localization and signal analysis}. IEEE Trans. Inform. Theory 36(5),  961--1005 (1990).




\bibitem{DLL} I. Daubechies,  H. J. Landau, Z. Landau, \textit{Gabor time–frequency lattices and the Wexler–Raz identity.}
J. Fourier Anal. Appl. 1(4), 437--478 (1995).





\bibitem{DSframe} R. J. Duffin,  A. C. Schaeffer: \textit{A class of nonharmonic Fourier
series}.   Trans. Am. Math. Soc.  72, 341--366 (1952).

\bibitem{Feicht83}
H. G. Feichtinger: \textit{A new family of functional spaces on the Euclidean n-space}. In: Proc.
Conf. on Theory of Approximation of Functions, Teor. Priblizh., 1983.

\bibitem{FG} H. G. Feichtinger, K. Gr\"{o}chenig: 
\textit{Gabor Frames and Time-Frequency Analysis
of Distributions}. J. Funct. Anal. 146, 464--495 (1997).


\bibitem{FK}
H. G. Feichtinger, W. Kozek: \textit{Quantization of TF lattice-invariant operators on elementary LCA
groups.}  In: Feichtinger H.G., Strohmer T. (eds) Gabor Analysis and Algorithms. Applied and Numerical Harmonic Analysis, pp. 233--266, Birkhäuser, Boston (1998).



\bibitem{FS1}
H. G. Feichtinger, T. Strohmer (Eds.): \textit{Gabor Analysis and Algorithms. Theory and Applications}.  
 Birkh\"{a}user, Basel (1998).



\bibitem{FS2} H. G. Feichtinger, T. Strohmer (Eds.): \textit{Advances in Gabor Analysis}.  
Birkh\"{a}user, Basel (2003).



\bibitem{Gbook} K. Gr\"{o}chenig: {\it Foundations of Time-Frequency Analysis}. 
Birkh\"{a}user, Boston (2000).


\bibitem{GrAlg} K. Gr\"{o}chenig:  \textit{Acceleration of the frame algorithm}. IEEE Trans. Signal Process. 41(12), 3331--3340 (1993).




\bibitem{HLS}
E. Hayashi, S. Li, T. Sorrells: \textit{Gabor duality characterizations}.
In: Heil C. (Ed.) Harmonic Analysis and Applications. Applied
and Numerical Harmonic Analysis, pp. 127--137. Birkh\"auser, Boston
(2006).
 


\bibitem{Hbook} C. Heil:  \textit{A Basis Theory Primer}. Expanded ed., Birkh\"auser, Basel, 2011.


\bibitem{Jannsen} A. J. E. M. Janssen:  \textit{Signal analytic proofs of two basic results on lattice expansions}. Appl. Comput. Harmonic Anal. 1, 350--354 (1994).

\bibitem{JanssenSchulzAlg} 
A. J. E. M. Janssen: \textit{Some iterative algorithms to compute canonical
windows for Gabor frames}. In: S. S. Goh, A. Ron, and Z. Shen (Eds.), Gabor and Wavelet Frames IMS Lecture Notes Series, vol. 10,  51--76 (2007).





\bibitem{JanssenCanDualWindow}
A. J. E. M. Janssen: \textit{Duality and biorthogonality for
Weyl–Heisenberg frames}. 
 J. Fourier Anal. Appl.   1(4), 403--436 (1994).

 

\bibitem{Tobias} 
T. Kloos, J. St\"{o}ckler, and K. Gr\"{o}chenig: 
\textit{Implementation of discretized Gabor
frames and their duals}. IEEE Transactions on Information Theory, Vol. 62, No. 5,  2759--2771 (2016).



\bibitem{Li} S. Li: \emph{On general frame decompositions}. Numer. Funct. Anal. Optim.
16(9-10), 1181--1191 (1995).





\bibitem{LLM} S. Li, Y. Liu, and T. Mi: \textit{Sparse dual frames and dual Gabor functions
of minimal time and frequency supports}.  
J. Fourier Anal. Appl. 19, 48--76 (2013). 



\bibitem{PHLB}
N. Perraudin, N. Holighaus, P. L. S\o{}ndergaard,  P. Balazs:
\textit{Designing Gabor windows using convex optimization}. 
Appl. Math. Comput. 330, 266--287 (2018).





\bibitem{Strohmer}
T. Strohmer: \textit{Approximation of dual Gabor frames, window
decay, and wireless communications}.  
Appl. Comput. Harmonic Anal. 11, 243--262 (2001).

\bibitem{WR} J. Wexler, S. Raz: \textit{Discrete Gabor expansions}.  
Signal Process. 21(3), 207--220 (1990). 

\end{thebibliography}
\end{document}